\documentclass{article}
\usepackage{graphicx}
\usepackage{amsfonts}

\usepackage{amsmath, amssymb, amsthm, dsfont, bm}
\newtheorem{theorem}{Theorem}
\newtheorem{definition}{Definition}

\begin{document}

\title{Nonlinear optimization in Hilbert space using Sobolev gradients with applications}

\author{P. Kazemi\thanks
        {Ulm, DE
         ({\tt parimah.kazemi@gmail.com})}, \and
         R. J. Renka\thanks
        {Department of Computer Science \& Engineering,
         University of North Texas,
         Denton, TX 76203-1366
         ({\tt renka@cs.unt.edu})}}

%\begin{document}
\maketitle

\begin{abstract}
The problem of finding roots or solutions of a nonlinear partial differential equation may be formulated as the problem of minimizing a sum of squared residuals.  One then defines an evolution equation so that in the asymptotic limit a minimizer, and often a solution of the PDE, is obtained.  The corresponding discretized nonlinear least squares problem is an often met problem in the field of numerical optimization, and thus there exist a wide variety of methods for solving such problems.  We review here Newton's method from nonlinear optimization both in a discrete and continuous setting and present results of a similar nature for the Levernberg-Marquardt method.  We apply these results to the Ginzburg-Landau model of superconductivity.
\end{abstract}

%\begin{keywords}
%  Gradient system, Levenberg-Marquardt, Sobolev gradient
%\end{keywords}

%\begin{AMS}
%   35A15
%\end{AMS}

\section{Introduction}   % SECTION 1
Consider the problem of finding a solution to the PDE
\begin{equation}\label{pde}
F(Du)=0
\end{equation}
where for $u$ in the Sobolev space $H:=H^{1,2}(\Omega)$, $Du = \{ D_{\alpha} u : |\alpha| \leq 1 \}$.  $\Omega$ is assumed to be a bounded domain  of dimension $n$ with smooth boundary.  $F$ is a function from $\mathbb{R}^{n+1}$ to $\mathbb{R}^m$ which is commonly referred to as a Nemistkii operator.  Let $L:=L^2(\Omega)$.  In order that $F \circ D :H \rightarrow [L]^m$ be Fr\'{e}chet differentiable it is sufficient that $F$ be $C^1$ and satisfy the growth bound $|F'(x)| \leq c |x|$ for $x \in \mathbb{R}^{n+1}$ and some $c>0$ (see \cite{chill}).  One can propose to find a solution of \eqref{pde} by solving the minimization problem
\begin{equation}\label{min}
E(u) = \frac{\|F(Du)\|_{L}^2}{2}. 
\end{equation}
The proposed method would ideally be an existence results, proving that a minimizer exists and is a solution of the PDE, and it should give a recipe for computing such a solution numerically as for most nonlinear problem closed form solutions are not possible.  In order to be effective, the numerical method should emulate an iteration in the infinite-dimensional Sobolev space in which the PDE is formulated. This formulation follows naturally from Neuberger's theory of Sobolev gradients.  The Fr\'{e}chet derivative $E'(u)$, which is a bounded linear functional on $H$,  is represented by an element of $H$.  This element is the Sobolev gradient of $E$ at $u$ 
and is denoted by $\nabla_H E(u)$:
\[  E'(u)h = \langle h,\nabla_H E(u) \rangle_H, \;\; h \in H.  \]
Note that the gradient depends on the inner product attached to $H$.
One considers the evolution equation
\begin{equation} \label{flow}
  z(0) = z_0 \in H \mbox{ and } z'(t)=-\nabla_H E(z(t)), \;\; t \geq 0.
\end{equation}
The energy $E$ is non-increasing on the trajectory $z$.  Existence,
uniqueness, and asymptotic convergence to a critical point are established
by the following two theorems taken from \cite[Chapter 4]{neuberger2010}.
\begin{theorem} \label{zexists}
Suppose that $E$ is a non-negative $C^1$ real-valued function on a
Hilbert space $H$ with a locally Lipschitz continuous Sobolev gradient.  
Then for each $z_0 \in H$ there is a unique global solution of \eqref{flow}.
\end{theorem}
\begin{definition}\label{gradineqdef}
The energy functional $E$ satisfies a \emph{gradient inequality} on 
$K \subseteq H$ if there exists $\theta \in (0,1)$ and $m > 0$ so that 
for all $x \in K$
\begin{equation*}
   \| \nabla_H E(x) \|_H \geq m E(x)^{\theta}.
\end{equation*}
\end{definition}
\begin{theorem}\label{gradineq}
Suppose that $E$ is a non-negative $C^1$ functional on $H$ with a
locally Lipschitz continuous gradient, $z$ is the unique global solution 
of \eqref{flow}, and $E$ satisfies a gradient inequality on the 
range of $z$.  Then $\lim_{t \rightarrow \infty} z(t)$ exists and is a 
zero of the gradient, where the limit is defined by the $H$-norm.  By
the gradient inequality, the limit is also a zero of $E$.
\end{theorem} 

The above theorems provide a firm theoretical basis for the numerical 
treatment of a system of nonlinear PDE's by a gradient descent method that
emulates \eqref{flow}; i.e., discretization in time and space results
in the method of steepest descent with a discretized Sobolev gradient.
Note that the Sobolev gradient method differs from methods based on
calculus of variations in which the Euler-Lagrange equation is solved.
Forming the Euler-Lagrange equation requires integration by parts to
obtain the element that represents $E'(u)$ in the $L^2$ inner product.
This $L^2$ gradient is usually only defined on a Sobolev space of higher
order than that of $H$.  Hence, unlike the Sobolev gradient, the $L^2$ 
gradient is only densely defined on the domain of $E$.  For gradient
flows involving the $L^2$ gradient, existence and uniqueness results
similar to those of Theorems \ref{zexists} and \eqref{gradineq} may
be proved, but only under stricter assumptions.

\section{Newton and Variable Newton methods}
In order to solve the minimization problem \eqref{min}, we take the first variation to obtain 
\begin{equation}\label{deriv}
E'(u)h=\langle F'(Du)Dh, F(Du) \rangle_L.
\end{equation}
We seek an evolution equation so that in the limit as time goes to infinity, we can find a zero.  A natural setting would be if such an evolution came from a gradient system as defined in \cite{chill}.  In particular, first assume that $G(u):= F(Du)$ and $G'(u)$ is invertible for each $u$.  Then Newton's method 
\begin{equation*}
u(0)=u_0 \text{ and } u'(t) = - (G'(u))^{-1}G(u)
\end{equation*}
is the gradient system associated with the inner product
\begin{equation}
\langle v , w \rangle_{g_N(u)} = \langle G'(u) v , G'(u) w \rangle_L. 
\end{equation}
This is achieved by noting that
\[ E'(u)h = \langle  G'(u) h , G'(u) (G'(u))^{-1} G(u) \rangle = \\
\langle h , (G'(u))^{-1} G(u) \rangle_{g_N(u)}. \]

Continuous Newton's method gives an infinite dimensional method for finding solutions of \eqref{pde} ; see, e.g., \cite{newton}, \cite{castro}, and \cite{neuberger2007a} for zero-finding results of Nash-Moser type (\cite{moser}).  In \cite{karatson07} Newton's method is discussed in relation to gradient descent methods.  It is shown that, while the method of steepest descent is locally optimal in terms of the descent direction for a fixed metric, Newton's method is optimal (in a sense which is made precise) in terms of both the direction and the inner product in a variable metric method. When Newton's method is available the quadratic rate of convergence to a solution makes this method ideal in a numerical setting. 

For many partial differential equations $G'(u)$ may not be invertible for all $u$ thus Newton's method cannot be applied in the infinite dimensional setting.  In the finite dimensional setting, one needs that the initial condition be close to the solution to obtain convergence. Another option is to minimize \eqref{min} using a variable metric method. We give here the description of one such method which when discretized gives a variation of the Levenberg-Marquardt method.  The results are taken from \cite{kazemirenkalm}.

For $u \in H=H^1(\Omega)$ consider the bilinear form on $H$ defined by
\begin{equation}\label{normu}
  \langle v,w \rangle_{u} = \langle v,w \rangle_{H} + (1/\lambda(u))
  \langle G'(u)v , G'(u)w \rangle_{L}
\end{equation}
where $\lambda(u)$ is a positive damping parameter.  By our assumption that $G'(u)  \in L(H,L)$, there exists a constant $c=c(u)$ so that for all $v \in H$, $\|G'(u)v\|_{L} \leq c \|v\|_{H}$.  Hence $\|v\|_u^2 =
\|v\|_H^2 + (1/\lambda(u))\|G'(u)v\|_{L}^2$, and
\begin{equation*}
  \|v\|_{H} \leq \|v\|_{u} \leq \sqrt{1 + c(u)^2/\lambda(u)}\; \|v\|_{H}
\end{equation*}
so that, for each $u \in H$, \eqref{normu} defines a norm that is equivalent to the standard Sobolev norm on $H$.  The gradient of $E$ with respect to $\langle \cdot , \cdot \rangle_u$ is defined to be the unique element $\nabla_u E(u)$ so that 
\begin{equation}\label{gradu}
E'(u)h  = \langle h , \nabla_u E(u) \rangle_{u} \ \forall \ h \in H.
\end{equation}
Consider the gradient flow
\begin{equation}\label{gradflow}
z(0) = u_0 \in H \mbox{ and } z'(t)=-\nabla_{z(t)} E(z(t)), \;\; t \geq 0.
\end{equation}
We seek a $C^1$ solution $z$ of this flow so that $u_f = \lim_{t \rightarrow \infty} z(t)$ exists and $E'(u_f) = 0$.  In the case that a gradient inequality is satisfied, a zero of the derivative of $E$ is also a zero of $E$ and hence a solution of \eqref{pde}.  The key idea in obtaining global existence and asymptotic convergence of the flow \eqref{gradflow} is to obtain an expression for the abstract gradient. We obtain this expression by considering a family of orthogonal projection onto the graph of a closed densely defined operator.  In particular, let $S_u : H \rightarrow [L]^{n+1}$ be given by $S_u=\binom{T} {T_u}$ where $Th=\{ D_{\alpha} h : |\alpha|=1 \}$ and $T_u h = (1/ \lambda(u)) G'(u)h$.  Since the domain of $S_u$ is all of $H$, $S_u$ can be viewed as a densely defined operator on $L$.  It is also the case that $S_u$ is a bounded 
linear operator from $H$ to $[L]^{n+1}$.  Note that $S_u$ need not be bounded when viewed as an operator from $L$ to $[L]^{n+1}$.  It also follows that the graph of $S_u$ is a closed subspace of $[L]^{n+2}$.  By a theorem of von Neumann (\cite{vonneumann}) there exists a unique orthogonal projection from $[L]^{n+2}$ onto the graph of $S_u$, and the projection is given by 
\begin{equation*}
   P_u = \left( \begin{array}{ll}
          (I+S_u^* S_u)^{-1}    &  S_u^*(I+S_u S_u^*)^{-1} \\
          S_u(I+S_u^* S_u)^{-1} &  I-(I+S_u S_u^*)^{-1}
          \end{array} \right).
\end{equation*}
This result can also be found in \cite[Theorem 5.2]{neuberger2010}.  Here $S_u^*$ is the adjoint of 
$S_u$ with $S_u$ treated as a closed and densely defined operator on $L$,
and hence $S_u^*$ is also closed and densely defined on its domain 
$[L]^{n+1}$.  Note also that $(I+S_u^* S_u)^{-1}$ and $S_u^*(I+S_u S_u^*)^{-1}$ are 
everywhere defined on $L$ and $[L]^{n+1}$, respectively, and are bounded 
as operators from $L$ to $L$ and $[L]^{n+1}$ to $L$, respectively 
(\cite[Sec 118]{rsznagy}).

\subsection{An expression for the gradient}  % SECTION 3.3

We will obtain an expression for the gradient given in \eqref{gradu}.
The graph of $S_u$ is $\left\{ \binom{Dh}
{T_uh} : h \in H \right\}$.  Since $P_u$ is the unique orthogonal 
projection of $[L]^{n+2}$ onto the graph of $S_u$, $P_u$ is the 
identity on the graph of $S_u$, and thus $P_u \binom{Dh}{T_uh} = 
\binom{Dh}{T_uh}$ for all $h \in H$.  Using symmetry of $P_u$, we have
\begin{eqnarray*}
  E'(u)h = \langle G'(u)h , G(u) \rangle_{L} =  \\ 
  \sqrt{\lambda(u)} \left\langle \binom{Dh}{(1/\sqrt{\lambda(u)})G'(u)h} , 
     \binom{0} {G(u)} \right\rangle_{[L]^{n+2}}=  \\
  \sqrt{\lambda(u)} \left\langle P_u \binom{Dh}{(1/\sqrt{\lambda(u)})G'(u)h}, 
     \binom{0} {G(u)} \right\rangle_{[L]^{n+2}}=  \\
  \sqrt{\lambda(u)} \left\langle \binom{Dh}{(1/\sqrt{\lambda(u)})G'(u)h} , 
     P_u \binom{0} {G(u)} \right\rangle_{[L]^{n+2}}=  \\
  \sqrt{\lambda(u)} \left\langle h , \Pi P_u \binom{0} {G(u)} \right\rangle_{u},
\end{eqnarray*}
where $\Pi$ is the operator that extracts the first element of a pair:
$\Pi \binom{x}{y} = x$.  For each $u \in H$ we have a gradient
\begin{equation}  \label{formula}
  \nabla_u E(u) = \sqrt{\lambda(u)} \Pi P_u \binom{0} {G(u)}= 
  \sqrt{\lambda(u)} S_u^*(I+S_uS_u^*)^{-1} \binom{0}{G(u)}.
\end{equation}
\begin{theorem}
\begin{equation} \label{graduexp}
  \nabla_{u} E(u) = \sqrt{\lambda(u)} M T_u^*( I + T_u M T_u^*)^{-1} F(Du),
\end{equation}
where $T_u^*$ is the adjoint of $T_u$ when viewed as a closed and densely defined operator 
on $L^2$, and $M = (I + T^* T)^{-1} = (D^* D)^{-1}$ is a smoothing operator.
\end{theorem}

A steepest descent iteration with a discretization of this gradient is a 
generalized Levenberg-Marquardt iteration in which the identity or a diagonal
matrix is replaced by the positive definite operator $M^{-1} = D^*D$.

The generalized Levenberg-Marquardt method is given by
\begin{equation*}
u_{n+1} = u_n - \lambda(u)(\lambda(u) D^*D +(G'(u))^* G'(u))^{-1}(G'(u))^*F(Du)
\end{equation*}
which is a forward Euler discretization of \eqref{gradflow} with time step 1.  The expression of the gradient was used to obtain the following results.

\begin{theorem}\label{exist}
Suppose $E$ is as defined in \eqref{min} with $F \circ D$ a $C^2$ function 
defined on $H$ with range in $L^2$, and suppose that $\lambda: H \rightarrow
\mathbb{R}$ is locally Lipschitz continuous and bounded below by a
positive constant.  Then the gradient system 
\eqref{gradflow} has a unique global solution $z \in C^1( [0,\infty) , H)$.
\end{theorem}
\begin{theorem}\label{pkrr}
Suppose that there exists $\xi \in (0,2)$ so that if $u \in H$, there is $\gamma(u) > \|F(Du)\|_{L}^{\xi}$ such that for each $g$ in the domain of $G'(u)^*$ with $\|g\|_L=1$ the linear PDE 
\[ G'(u) x = g \]
has a solution $x \in H$ with $\|x\|_H^2 \leq \frac{1}{\gamma(u)}$.  Then a gradient inequality is satisfied.  Here $G'(u)^*$ denotes the adjoint of $G'(u)$ when viewed as a densely define operator on $L$.
\end{theorem}

\begin{theorem}\label{main}
Suppose that the hypotheses of Theorem \ref{exist} are satisfied 
so that \eqref{gradflow} has a unique global solution $z$, and
suppose that the hypotheses of Theorem \ref{pkrr} are satisfied 
on an open region containing the range of $z$.  Then 
$u = lim_{t \rightarrow \infty} z(t)$ exists and $F(Du) = 0$.
\end{theorem}

\section{Applications to superconductivity}
\begin{figure}[h!]
\centering
\includegraphics[width=0.45\columnwidth]{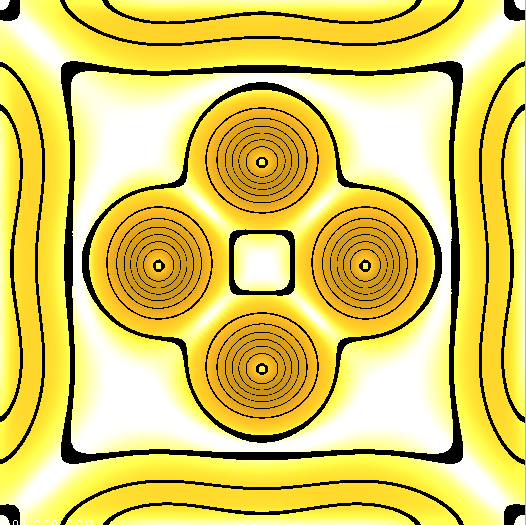} 
\includegraphics[width=0.45\columnwidth]{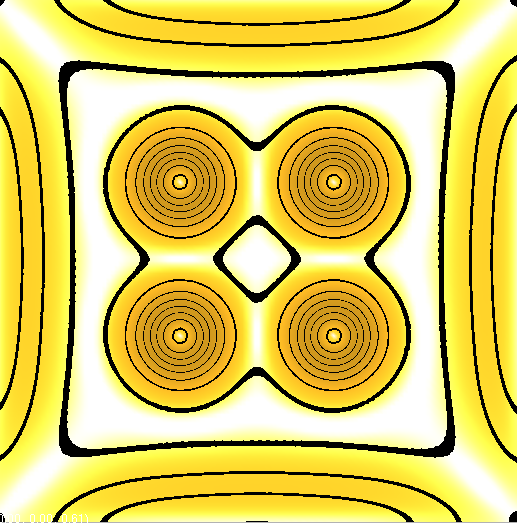}\\
\includegraphics[width=0.45\columnwidth]{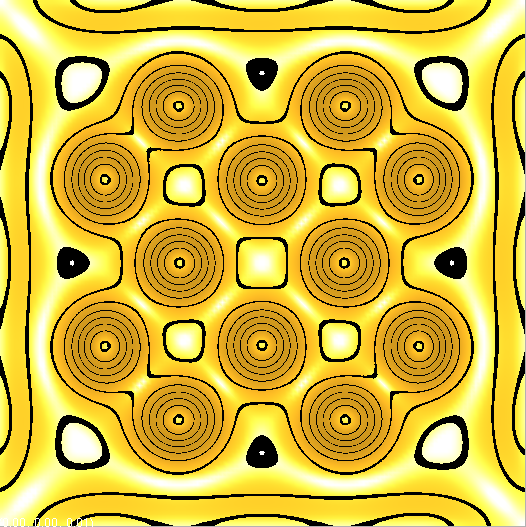} 
\includegraphics[width=0.45\columnwidth]{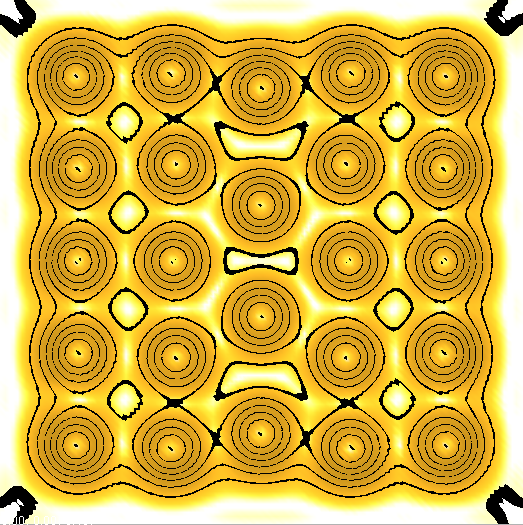}
\caption{Vortex configurations corresponding to a density plot of the minimizer for $\kappa=4$ and $H_0=4,4,6,8$.}\label{figure} 
\end{figure}

In \cite{kazemirenkagl}, we applied a variation of the above formulation to study the Ginzburg-Landau energy.  The Ginzburg-Landau model postulates that the behavior of the superconducting electrons in materials can be described by a complex valued wave function in which case the Ginzburg-Landau (Gibbs free) energy is given by
\begin{equation}\label{eqngl}
E(u,A)= \int_{\Omega} \frac{|\nabla u - i Au |^2}{2} + \frac{| \nabla \times A - H_0 |^2}{2} + \frac{\kappa^2}{4} ( |u|^2 - 1)^2 
\end{equation}
in nondimensionalized form.  Here $u$ is the complex valued wave function which gives the probability density of the superconducting electrons, $A$ is the induced magnetic vector potential, $H_0$ is the external magnetic field, and $\kappa$ is the Ginzburg-Landau coefficient which characterizes the type of superconducting sample (type I or type II). The central Ginzburg-Landau problem is to find a minimizer of the Ginzburg-Landau energy.  Note that this energy can be written in the form \eqref{min} with

\begin{equation}
F(D(u,A))=
\begin{pmatrix}
r_1 + as \\
s_1 - a r\\
r_2 + b s\\
s_2 - b r\\
b_1 - a_2 - H_0\\
\frac{\kappa}{\sqrt{2}} (r^2 +s^2 -1)
\end{pmatrix}.
\end{equation}
for $u=\binom{r}{s}$ and $A=\binom{a}{b}$.  One can check that $F \circ D$ is $C^2$ from $H= [H^1(\Omega)]^4$ to $L=[L^2(\Omega)]^6$.  Thus the gradient system \eqref{gradflow} has a unique global solution when $\lambda$ satisfies the properties of Theorem \ref{exist}.  We studied the rate of convergence of this flow to a minimizer using a trust-region method in \cite{kazemirenkagl}.  In future work, it would be a very nice result to obtain a gradient inequality to prove convergence by verifying the condition of Theorem \ref{pkrr}.  Since for the Ginzburg-Landau energy, a minimizer does not correspond to a zero of the energy, the definition of the gradient inequality is altered to the following definition taken from \cite{chill}.
\begin{definition}
Suppose $E$ is as in \eqref{min} and that $E$ achieves a local minimum at $u_m$.  Then $E$ is said to satisfy a gradient inequality in a neighborhood of $u_m$ if there exists a ball $B$ containing $u_m$ and $\xi \in (0,1)$, $c>0$ such that all $v \in B$
\begin{equation*}
|E(v) - E(u_m)|^{\xi} \leq c \| \nabla_v E(v) \|_v.
\end{equation*}
\end{definition}
This formulation was used to obtain a stabilization result for the Ginzburg-Landau equations in \cite{du} and \cite{takac}.  In figure \ref{figure}, we give contour plots of minimizers for various parameters.

\section{Conclusion}\label{conclusion}
 
We extended the theory of Sobolev gradients to include gradients
associated with a variable inner product, and we described a generalized
Levenberg-Marquardt method as a gradient flow in an infinite-dimensional 
Sobolev space.  We presented conditions under which the flow is
guaranteed to converge to a zero of a residual representing a solution of 
a nonlinear partial differential equation.  The conditions include 
smoothness of the residual and satisfaction of a gradient inequality.

Our results provide a theoretical basis for a practical and effective
method of solving nonlinear partial differential equations.
As a further development we will seek to apply our results to particular
problems.  A proof that a numerical iteration has a convergent 
counterpart in the infinite-dimensional setting represents an important
contribution to the goal of a unified approach to treating partial
differential equations by numerical and analytical methods.

\section*{Acknowledgements}
The first author would like to thanks the Institute of Applied Analysis and the Institute of Quantum Physics at Ulm University for their hospitality during my time in Ulm.


\begin{thebibliography}{99}
%%
\bibitem{castro} A. Castro and J. W. Neuberger, An inverse function theorem 
  by means of continuous Newton's method, J. Nonlinear Analysis, 47 (2001), 
  pp. 3259--3270.
%%
\bibitem{chill} Ralph Chill, Eva Fasangova, Gradient Systems, 13th 
  International Internet Seminar, http://isem.univ-metz.fr.
%%
\bibitem{du} Fang-Hua Lin and Qiang Du, Ginzburg-Landau Vortices: Dynamics, Pinning, and Hysteresis, SIAM J. Math. Anal., 28 (6) (1997), pp. 1265--1293.
%%
\bibitem{newton} L.V. Kantorovich and G.P. Akilov, Functional Analysis, 
  Pergamon Press, Oxford, UK, 1982.
%%
\bibitem{karatson07} J Kar\'{a}tson and J. W. Neuberger, Newton's method 
  in the context of gradients,  Elec. J. Differential Equations 2007 (2007), 
  pp. 1--13.
%
\bibitem{kazemirenkagl} P. Kazemi and R. J. Renka, Minimization of the 
  {G}inzburg-{L}andau energy functional by a {S}obolev gradient trust-region 
  method, submitted.
%
\bibitem{kazemirenkalm} P. Kazemi and R. J. Renka, A Levenberg-Marquardt method based on Sobolev gradients, in preparation.
% 
\bibitem{moser} J. Moser, A rapidly convergent iteration method and 
  nonlinear differential equations, Ann. Scuola Normal Sup. Pisa, 20 
  (1966), pp. 265--315.
%  
\bibitem{neuberger2010} Neuberger J. W., Sobolev Gradients and Differential 
  Equations, 2nd edition, Springer, 2010.
%%
\bibitem{neuberger2007a} J. W. Neuberger, The continuous {N}ewton's method, 
  inverse functions and {N}ash-{M}oser, American Math. Monthly, 114 (2007), 
  pp. 432--437.
%
\bibitem{rsznagy} Riesz and Sz. Nagy, Functional Analysis, Dover Publications, Inc., 1990.
%  
\bibitem{takac} Eduard Feireisl and Peter Tak\'{a}c, Long-Time Stabilization of Solutions to the Ginzburg–Landau Equations of Superconductivity, Monatsh. Math., 133 (2001), pp. 197–221.
%
\bibitem{vonneumann} J. von Neumann, Functional Operators II, Annls. Math. 
  Stud., 22, 1940.
%
\end{thebibliography}
\end{document}